%% file: kneser-compare-jdmv-stripped.tex
\title{Topological lower bounds for the chromatic number: A~hierarchy
}
\date{\small November 12, 2003}

\documentclass[11pt]{article}

\newcommand{\cmt}[1]{\ifhmode\newline\fi{\sf *** \ \ #1 \\}}

\usepackage{a4wide,url,amsmath,amssymb,amsthm}
\usepackage{latexsym}
\usepackage{graphicx,color,paralist}
\newtheorem{theorem}{Theorem}

\newtheorem{proposition}[theorem]{Proposition}

\theoremstyle{definition}
\newtheorem{definition}[theorem]{Definition}
\theoremstyle{remark}

\newcommand\ie{i.\,e.}
\newcommand\eg{e.\,g.}
\def\proofend{\unskip \nobreak \hskip0pt plus 1fill \qquad$\Box$\medskip\noindent}

\def\proof{\par\medbreak\noindent{\bf Proof.}\enspace}

\def\lsoft{{l\kern-0.035cm\char39\kern-0.03truecm}}
\newcommand\dolnikov{Do{\lsoft}nikov}
\newcommand\kriz{K\v{r}\'{\i}\v{z}}
\newcommand\lju{Lyusternik}
\newcommand\snirel{Shnire{\lsoft}man}

\newcommand{\FF}{{\cal F}}
\newcommand{\CC}{{\cal C}}
\renewcommand{\AA}{{\cal A}}
\newcommand{\BB}{{\cal B}}

\newcommand{\K}{{\sf K}}
\newcommand{\B}{{\sf B}}
\renewcommand{\L}{{\sf L}}
\newcommand{\N}{{\sf N}}

\newcommand{\Boxc}[1]{\B(#1)}
\newcommand{\Bzero}[1]{\B_0(#1)}
\newcommand{\Bedge}[1]{\B_{\rm edge}(#1)}
\newcommand{\Bone}[1]{\B_{\rm chain}(#1)}
\newcommand{\Boxckneser}[1]{\B^{\rm{KG}}(#1)}
\newcommand{\Bzerokneser}[1]{\B_0^{\rm{KG}}(#1)}
\newcommand{\Bsark}[1]{\B^{\rm{KG}}_{\rm Sark}(#1)}
\newcommand{\Bonekneser}[1]{\B_{\rm chain}^{\rm{KG}}(#1)}

\newcommand{\Nbhd}[1]{\N(#1)}
\newcommand{\Lov}[1]{\L(#1)}
\newcommand{\Ebar}{\overline E}
\newcommand{\ebar}{\overline e}

\newcommand\ind{{\mathop {\rm ind\,}\nolimits}}
\newcommand\conn{{\mathop {\rm connectivity}\nolimits}}
\newcommand\acyc{{\mathop {\Z_2\mbox{\rm-acyclicity}}\nolimits}}
\newcommand\sd{{\mathop {\rm sd\,}\nolimits}}
\def\id{\hbox{\rm id}}
\newcommand\conv{{\mathop {\rm conv}\nolimits}}

\newcommand\susp{{\mathop {\rm susp\,}\nolimits}}
\newcommand\cd{{\mathop {\rm cd}\nolimits}}
\newcommand\KG{{\mathop {\rm KG}}}
\newcommand\SG{{\mathop {\rm SG}}}
\newcommand\CN{{\mathop {\rm CN}}}

\newcommand\stab{\mathrm{stab}}

\newcommand\sep{:\,}
\def\:{\colon}

\newcommand\assign{:=}
\newcommand\polyh[1]{{\|#1\|}}

\newcounter{eqncount}  
\def\theeqncount{{\rm H\arabic{eqncount}}}
\newcommand\eqntag[1]{\limits%
^{\refstepcounter{eqncount}\label{#1}\rm(\theeqncount)}}
\newcommand{\R}{{\mathbb R}}
\newcommand{\Z}{{\mathbb Z}}

\newcommand{\heading}[1]{\vspace{1ex}\par\noindent{\bf #1}}

\author{
 {\sc Ji\v{r}\'{\i} Matou\v{s}ek}\\
{\footnotesize Department of Applied Mathematics and}\\[-1.5mm]
{\footnotesize Institute for Theoretical Computer Science (ITI)}\\[-1.5mm]
{\footnotesize  Charles University, Malostransk\'{e} n\'{a}m. 25}\\[-1.5mm]
{\footnotesize  118~00~~Praha~1, Czech Republic}\\[-1.5mm]
{\footnotesize \url{matousek@kam.mff.cuni.cz}}
\and
{\sc G\"unter M. Ziegler}\\
{\footnotesize Institute of Mathematics, MA 6-2}\\[-1.5mm]
{\footnotesize Technical University Berlin}\\[-1.5mm]
{\footnotesize D-10623 Berlin, Germany}\\[-1.5mm]
{\footnotesize \url{ziegler@math.tu-berlin.de}}
}

\begin{document}

\maketitle

\begin{abstract}
This paper is a study of ``topological''
lower bounds for the chromatic number of a graph.
Such a lower bound was first introduced by
Lov\'asz in 1978, in his famous proof of the
\emph{Kneser conjecture} via Algebraic Topology.
This conjecture stated that the
\emph{Kneser graph} $\KG_{m,n}$, the graph with all
$k$-element subsets of $\{1,2,\ldots,n\}$ as vertices and all pairs of
disjoint sets as edges, has chromatic number $n-2k+2$. 
Several other
proofs have since been published (by B\'ar\'any, Schrijver, \dolnikov,
Sarkaria, \kriz, Greene, and others), all of them based on 
some version of the
Borsuk--Ulam theorem, but otherwise quite
different.  Each can be extended to yield some lower bound on the
chromatic number of an arbitrary graph.
(Indeed, we observe that \emph{every} finite graph may be
represented as a generalized Kneser graph, to which the above
bounds apply.)

We show that these bounds are almost linearly
ordered by strength, the strongest one being essentially Lov\'asz'
original bound in terms of a  neighborhood complex.
We also present and compare
various definitions of a \emph{box complex} of a graph
(developing ideas of Alon, Frankl, and Lov\'asz and of \kriz).
A suitable box complex is equivalent to Lov\'asz'
complex, but the construction is simpler and functorial, mapping
graphs with homomorphisms to $\Z_2$-spaces with $\Z_2$-maps.
\end{abstract}

\section{Introduction}

Graph coloring is a classical combinatorial topic:
For a given (finite) graph $G$, determine how to distribute
a minimal number of colors to the vertices in such a way
that adjacent vertices get different colors. The minimum
number of colors is $\chi(G)$,
the \emph{chromatic number} of the graph.
The graph coloring
problem has numerous important practical motivations;
among the more recent ones, we mention that it appears as a
(simplified) model for the frequency assignment problem
in mobile communication
(cf.\ Bornd\"orfer et.\,al.\ \cite{borndoerfer98:_frequen}
and Eisenbl\"atter et al.\ \cite{eisenblaetter02:_frequen}).

The most famous graph coloring problem is, of course,
the Four Color Problem, asking whether every planar graph
can be colored by four colors, which was answered positively
by Haken and Appell 1977 and re-solved by Robertson, 
Sanders, Seymour \& Thomas
\cite{robertson97-FCT}. Even for planar graphs, though, 
determining 3-colorability is already algorithmically difficult
(NP-hard), and   beyond the range of planar graphs, 
the gaps between the upper and the lower bounds that one can 
reasonably obtain for $\chi(G)$ may be huge.
This claim is supported both by theoretical results 
on the complexity of graph coloring algorithms (see 
the discussion in Khanna, Linial \& Safra \cite{khanna00})
and from the perspective of combinatorial
optimization (see Mehrotra \& Trick \cite{MTrick}).

For any given graph $G$, we obtain an 
\emph{upper bound} on the chromatic number $\chi(G)$
by ``guessing'' a coloring, for example,
by  running a coloring heuristic.
(See \cite{culberson:_graph_color_page} for implementations,
and Reed \& Molloy \cite{molloy01:_graph_color_probab_method}
for a theoretical study of randomized algorithms.)
How can we check that such an upper bound is good,
that is, close to the actual chromatic number?
We need a \emph{lower bound} on the chromatic number.

Two other basic graph parameters lead to straightforward
``combinatorial''
lower bounds for $\chi(G)$: The \emph{clique number}
$\omega(G)$, the largest number of mutually adjacent
vertices in $G$, obviously
satisfies $\omega(G)\leq \chi(G)$, and 
we also have $\chi(G)\geq n/\alpha(G)$ for every graph $G$
on $n$ vertices, where $\alpha(G)$ is the \emph{independence
number} of $G$, that is, the maximum number of mutually nonadjacent
vertices in $G$.
Both $\omega(G)$ and $\alpha(G)$ are hard to compute, or even
approximate, for general graphs.
But even leaving this aside,  both of these bounds may be
very weak, since one can construct
graphs where both $\omega(G)$ and $n/\alpha(G)$
are arbitrarily small compared to  $\chi(G)$.

A considerably more sophisticated lower bound for $\chi(G)$
(still hard to compute in general)
is the \emph{fractional chromatic number} $\chi_f(G)$,
which can be briefly defined as the
minimum  ratio $a/b$ such that $G$ has a $b$-fold covering
by $a$ independent sets. (The name is motivated by
a slightly different but equivalent definition,
in which each color can still be used only on an independent
set of vertices, but one is allowed to use fractional
amounts of colors, say to color a vertex 
by $\frac 15$ of red and by $\frac 45$ of blue.)
The gap between $\chi_f(G)$ and $\chi(G)$
can still be arbitrarily large, but it is much harder
to come up with examples of this.
(See \cite[Part~VIII]{molloy01:_graph_color_probab_method}
for further discussion.)

Actually, only very few types of such examples
are known, and the arguably most important ones are
provided by \emph{Kneser graphs}.
The Kneser graph $\KG_{n,k}$ has all the $N=\binom nk$
 $k$-element subsets of $\{1,2,\ldots,n\}$ 
as vertices and all pairs of disjoint sets as edges.
It arose in an innocent little
problem that Martin Kneser posed in 1955 {here},
in the \textsl{Jahresbericht der DMV} \cite{Kneser}.
(Apparently the problem arose from Kneser's study of
a number-theoretic paper \cite{kaplansky53:_quadr} by Kaplansky.) 
Kneser asked for a proof  that
$\chi(\KG_{n,k})\geq n-2k+\mathit{const.}$,
and conjectured that $\chi(\KG_{n,k})= n-2k+2$, where 
$\chi(\KG_{n,k})\leq n-2k+2$ is easy to verify using
a simple greedy coloring.

Kneser's question posed a substantial challenge since
all the classical lower bounds listed above
fail for Kneser graphs. For suitable
parameters, say for $n=3k-1$, we have:
\begin{compactitem}[$\bullet$]
\item
The chromatic number is large,
$\chi(\KG_{3k-1,k})=k+1$, and an optimal coloring is
easy to find (by a greedy approach).
\item
At the same time, the clique number is small, $\omega(\KG_{3k-1,k})=2$
(the graph is triangle-free for $n<3k$).
\item
The independence number is huge,
$\alpha(\KG_{3k-1,k})=\binom{n-1}{k-1}$,
and thus the corresponding lower bound for the chromatic
number, which also happens to agree with the fractional chromatic
number, are small:
$N/\alpha(\KG_{3k-1,k})=\chi_f(\KG_{3k-1,k})=\frac{3k-1}k<3$.
\end{compactitem}
Other, more ``algebraic'' types of lower bounds on
the chromatic numbers of graphs, in terms of 
the Lov\'asz theta function 
(see Lov\'asz \cite{Lov} and Knuth \cite{knuth94_sandwich})
or on the eigenvalues of the adjacency matrix
(see van Lint \& Wilson \cite[Chap.~31]{lint92:_cours_combin} and
Godsil \& Royle \cite[Chap.~9]{godsil01:_algeb_graph_theor})
do not help to close the gap in the case of Kneser graphs.

In 1978 Lov\'asz \cite{Lovasz-kneser}
settled Kneser's conjecture by an
original application of a tool from Algebraic Topology,
the Borsuk--Ulam theorem. He thus  provided
a completely new type of lower bound which
in the case of the Kneser conjecture was tight.

In subsequent 
years a number of new proofs of the
Kneser conjecture and of various extensions of it became available.
All of them are of topological nature, and
all of them depend on the Borsuk--Ulam theorem or
some variant of it. (This remains essentially true 
despite recent demonstrations
that some proofs  can be combinatorialized
 \cite{Mat1}, \cite{Z77}, since
the underlying ideas are still topological.)
The proof methods are  diverse, though, and at first sight
they look almost unrelated.

The Kneser graphs $\KG_{n,k}$ are quite special, but
the methods known for proving Kneser's conjecture, starting with
Lov\'asz's \cite{Lovasz-kneser} break-through,
extend beyond the original examples:
Each of them yields, explicitly or implicitly, a lower bound
for the chromatic number of any graph
(as we will see, \emph{every} graph is a ``generalized Kneser graph''),
although they
are rather weak for some classes of graphs.
These bounds are, of course,
all tight in the case of the Kneser graphs $\KG_{n,k}$,
but otherwise their strengths are apparently different:
For example, only some of them yield tight bounds for
the chromatic number of certain subgraphs of the Kneser
graphs investigated by Schrijver \cite{Schrijver-kneser}.

In the following, we will compare the
various ``topological lower bounds on the chromatic number.''
In the course of our arguments, we also sketch
proofs for them; some of these
proofs are quite elementary and simple, 
and sometimes considerably simpler than the original
derivations.

We show that, surprisingly, the topological lower bounds
for the chromatic number resulting from  known proofs
fall neatly into a  hierarchy, which is essentially linearly ordered.
So, we show that
``(the index version of) the Lov\'asz' bound is stronger than the
(generalized) Sarkaria bound, which is stronger than the (generalized)
 B\'ar\'any bound, and also stronger than
the \dolnikov--\kriz\ bound.'' 
(This, of course, should not indicate
any comparison of the usefulness of the methods or of the interest
of the papers. The hierarchy concerns only the
bounds that, \emph{in principle},
can be obtained by straightforward generalizations of these approaches,
and it does not say anything about the feasibility of actually
obtaining these bounds.)

\heading{Lov\'asz' proof and box complexes. }
A slightly modernized version of Lov\'asz' original proof
works along the following lines. To every graph $G$,
one assigns a topological space $T(G)$. The construction
goes via a simplicial complex. This,
on the one hand, is a purely combinatorial object (a hereditary
set system), and on the other hand, it is canonically associated
with a topological space (the geometric realization).
Then, for every coloring of $G$ by $m$ colors, one constructs
a continuous map from the space $T(G)$ to
the space $T(K_m)$ assigned to the complete graph $K_m$.
To show that $G$ has no $m$-coloring, it suffices
to exclude the existence of a continuous map $T(G)\to T(K_m)$.

Of course, things cannot be as simple as this, since every
topological space has a continuous map to any nonempty
topological space, namely, a constant map. One has to 
consider extra structure on the space $T(G)$, called 
a \emph{$\Z_2$-action}. A $\Z_2$-action
on a topological space $T$ is a homeomorphism $\nu\:T\to T$
such that $\nu(\nu(x))=x$ for every $x\in T$. A primary
example of a topological space with a $\Z_2$-action,
or a \emph{$\Z_2$-space}, is the $n$-dimensional unit sphere
$S^n$ in $\R^{n+1}$ with the $\Z_2$-action given by
$x\mapsto -x$, i.\,e., the antipodality.

The cleverly constructed space $T(G)$ comes equipped with
a $\Z_2$-action $\nu=\nu_G$, and 
for $G=K_m$, the space $T(K_m)$ even miraculously
happens to be (equivalent to) a sphere with the antipodality
as the $\Z_2$-action! 
(The dimension of this sphere depends on $m$, and as we
will see, it may differ slightly in various possible
constructions of $T(G)$.)
Moreover, the continuous map $f\: T(G)\to T(K_m)$
obtained from an $m$-coloring of $G$ is a \emph{$\Z_2$-map},
meaning that it commutes with the $\Z_2$-actions:
$f(\nu_G(x))=\nu_{K_m}(f(x))$.

Here the Borsuk--Ulam theorem enters. The most popular version
of it states that for every continuous map $f\:S^n\to \R^n$,
there is a point $x\in S^n$ with $f(x)=f(-x)$. However, an
equivalent version asserts that there is no $\Z_2$-map
$S^n\to S^m$ for $m<n$. A suitable generalization of
this result can sometimes be used to show that, for a particular
graph $G$, there is no $\Z_2$-map of $T(G)$ into $T(K_m)$,
and consequently, that $G$ is not $m$-colorable.

This high-level outline of Lov\'asz' 
does not tell one how to construct
suitable spaces $T(G)$. We should also
remark that Lov\'asz' original proof proceeded in a slightly
different way, without explicitly introducing the $\Z_2$-map.
A class of constructions which allows to phrase the proof
in the simple and conceptual way sketched above are
various \emph{box complexes} assigned to a graph,
which can successfully
play the role of $T(G)$ in the above discussion.
A box complex first appears in Alon, Frankl, and 
Lov\'asz \cite{AlonFranklLovasz} (where it was defined for hypergraphs),
and another version was used by 
\kriz~\cite{Kriz92}.

Formal definitions of box complexes will be given later, but roughly
we can say that a box complex of a graph $G$ is made of all complete
bipartite subgraphs of $G$.  A complete bipartite subgraph of $G$ is
specified by two disjoint subsets $A'$ and $A''$ of the vertex set,
such that every vertex $a'\in A'$ is connected by an edge to every
vertex $a''\in A''$.  The exchange of the subsets $A'$ and $A''$
yields an involution on the box complex, and makes it into a
$\Z_2$-space.

On an intuitive level,
the construction of the continuous $\Z_2$-map from the
box complex of $G$ into the box complex of $K_m$ can be described
as follows (we are indebted by  Lov\'asz for a beautiful
summary of the proof, which inspired much of the
present discussion): If $c$ is a proper $m$-coloring of $G$,
then whenever two disjoint sets $A',A''$ determine a complete
bipartite subgraph in $G$, they are assigned two disjoint
color sets $c(A')$ and $c(A'')$, which thus determine
a complete bipartite subgraph of $K_m$.
Furthermore, if $B'\supseteq A'$ and $B''\supseteq A''$
give a larger complete bipartite subgraph,
they receive larger color sets $c(B')\supseteq c(A')$ and
$c(B'')\supseteq c(A'')$. This allows one to define the
continuous map of the box complexes. Finally,
if we interchange $A'$ and $A''$, the color sets are
interchanged as well, and this makes the map of the box complexes
a $\Z_2$-map. 

Summarizing, an $m$-coloring of $G$ yields a $\Z_2$-map of the
box complex of $G$ into the box complex of $K_m$. If we use
the convenient notion of the \emph{index} of a $\Z_2$-space,
which is the smallest $m$ such that the $\Z_2$-space can
be $\Z_2$-mapped into the sphere $S^m$ with the antipodal
$\Z_2$-action, we get that the index of the box complex
of any $m$-colorable graph has to be at least as large as
the index of the box complex of $K_m$. The latter can be computed
once and for all (as we remarked above, the box complex
of $K_m$ happens to be equivalent to a sphere with the antipodal
action). Thus, the application of this method boils down
to bounding below the index of the box complex of $G$.

Interestingly, these ideas have several different implementations:
there are several distinct possibilities to define ``box complexes.''
The different box complexes have different ground sets, they
are of different sizes (the numbers of vertices/faces
differ on an exponential scale!), and some of them
may be considerably easier to use than others.

Our current favorite is the box complex $\Boxc{G}$, defined
in Section~\ref{s:overw} below. It has a small vertex set
(the disjoint union of two copies of $V(G)$),
and it yields the strongest bounds available.
However, we invite the reader to survey the panorama and to make
his/her own choices---several more versions of box complexes 
are discussed in Section~\ref{sec:box_complexes}.
We will show that many of them are equivalent for the
purposes of estimating the chromatic number.

\section{Preliminaries}\label{s:prelim}

Here we recall some general notions, facts, and notation
needed for a precise
statement of the results. We repeat some of the
definitions mentioned in the introduction more formally.

Unfortunately, for space reasons, we cannot afford to
introduce all the required topological notions at a leisurely
pace, and so the rest of the paper may not be easily accessible
without some knowledge of Topological Combinatorics.
We refer to \cite{Bjorner-topmeth} or \cite{Zivaljevic-topmeth}
for surveys of the terminology
and tools employed in this paper and to
 \cite{Matousek:BU} for a detailed textbook treatment.
On the other hand, 
readers with basic knowledge of the area may perhaps want
to skip this section and refer to it as needed during further reading.

\heading{Graphs. } The vertex set of a graph $G$ is written as $V(G)$,
and the edge set as $E(G)$. We suppose that all graphs are finite, simple,
and undirected. In order to avoid some trivial special cases,
 we also assume that the considered graphs have {no isolated vertices}.

 A \emph{homomorphism} of a graph $G$
into a graph $H$ is a mapping $f\:V(G)\to V(H)$ that preserves
edges; that is, $\{f(u),f(v)\}\in E(H)$ whenever $\{u,v\}\in E(G)$.
For our purposes, it is  convenient to regard a (proper) coloring
of $G$ by $m$ colors as a homomorphism of $G$ into the complete graph
$K_m$. The chromatic number of $G$ is denoted by $\chi(G)$.

We regard a \emph{bipartite graph} as a triple $(V',V'',E)$,
where $V',V''\subseteq V$ are disjoint and
$E\subseteq \{\{v',v''\}\sep v'\in V', v''\in V''\}$.
For a bipartite graph we assume
that the bipartition is fixed, and we also 
distinguish $(V',V'',E)$ from $(V'',V',E)$.
Since we will be concerned with various colorings of graphs,
we call  the parts $V'$ and $V''$
of the bipartition  the \emph{shores},
rather than the more common ``color classes.''
If $A',A''$ are disjoint subsets of the vertex set of some graph
$G$, we write $G[A',A'']$ for the \emph{bipartite} subgraph
with shores $A'$ and $A''$ induced by~$G$. (Note that this is not
necessarily an induced subgraph of~$G$, since only edges between
distinct shores are included.)

\heading{Kneser graphs. }
Let $X$ be a finite set and $\FF\subseteq 2^X$ a system of subsets of~$X$.
The \emph{Kneser graph} $\KG(\FF)$  has vertex set
$\FF$, and the edges  are all pairs of disjoint sets in~$\FF$.
For notational convenience, we assume that $X=[n]\assign \{1,2,\ldots,n\}$,
unless stated otherwise. Kneser's conjecture can be succinctly stated as
$\chi(\KG(\binom{[n]}{k}))= n-2k+2$ for $n\geq 2k> 0$,
where $\binom{[n]}k$ denotes the family of all $k$-element
subsets of $[n]$. In particular, $\KG_{n,k}=\KG(\binom{[n]}{k}))$.

It is easy to see that \emph{every} (finite) graph $G=(V,E)$ 
can be represented
as a Kneser graph of some set system. A simple and natural representation
is this: Let $\Ebar:=\binom{V}2\setminus E$ denote the set of non-edges of~$G$,
and for every $v\in V$, let us set $F_v\assign\{\ebar\in\Ebar\sep v\in\ebar\}$.
The Kneser graph of $\{F_v\sep v\in V\}$ is isomorphic to $G$;
the only problem is that the sets $F_v$ need not be all distinct
(for example, for $G=K_n$, we have $F_v=\emptyset$ for all $v$).
To remedy this, one can define $F'_v\assign F_v\cup \{v\}$,
obtaining distinct sets.
For a more economical representation, we can let $\CC$ be a covering
of $\Ebar$ by cliques (each $C\in\CC$ is a complete subgraph
of $(V,\Ebar)$ and each edge of $\Ebar$ is contained in some
$C\in\CC$). For $v\in V$, we then define 
$F''_v\assign\{C\in\CC\sep v\in C\}$; this is a potentially much
smaller Kneser representation. The problem of finding
a Kneser representation with the smallest ground set,
\ie, the smallest $\CC$, is the minimum clique cover for the
complement of $G$, and hence NP-complete and 
hard to approximate;
see, \eg,  Ausiello et al.\ \cite{ausiello99:_compl_approx}.

\heading{Simplicial complexes. }
We use letters like $\K,\L,\ldots$ to denote simplicial complexes.
(See, \eg, \cite{Munkres}, \cite{Bjorner-topmeth}, \cite{Matousek:BU} for
more background).
We consider only finite simplicial complexes, so
a simplicial complex $\K$ is a nonempty hereditary set system
(\ie, $S\in\K$ and $S'\subset S$ implies $S'\in\K$);
in particular, $\emptyset\in\K$. 
For example, $2^{[n]}$ is the $(n-1)$-dimensional simplex
considered as a simplicial complex.
We let $V(\K)$ denote the vertex set of $\K$, and
$\polyh{K}$ denotes the polyhedron of $\K$
(but sometimes we write just $\K$ for the polyhedron too,
when it is clear that we mean a topological space).
The \emph{dimension} of the complex $\K$ 
is $\dim\K\assign\max\{|S|{-}1\sep S\in\K\}$.
A \emph{simplicial map} of a simplicial complex $\K$
to a simplicial complex $\L$ is a map $f\:V(\K)\to V(\L)$  
such that $f(S)\in\L$  for all $S\in\K$.

For a partially ordered set $(X,\preceq)$,
the \emph{order complex} $\Delta(X,\preceq)$
has $X$ as the vertex set and all chains as simplices;
that is, a simplex has the form $\{x_1,x_2,\ldots,x_k\}\subseteq X$
with $x_1\prec x_2\prec\cdots\prec x_k$.
In particular, if $\FF$ is a set system, we write
$\Delta\FF$ for $\Delta(\FF\setminus\{\emptyset\},\subseteq)$.
If $\K$ is a simplicial complex, then 
$\Delta\K$ is the \emph{first barycentric subdivision} of $\K$,
also denoted by $\sd \K$ (the empty simplex $\emptyset$ is not a vertex
of the barycentric subdivision, and this is the reason for
removing $\emptyset$ in the definition of $\Delta\FF$).

The  (twofold) \emph{deleted join} of $\K$,
denoted by $\K^{*2}_\Delta$, has vertex set
$V(K){\times} [2]$ (two copies of $V(\K)$) and
the simplices are $\{ S_1\uplus S_2\sep S_1,S_2\in\K, S_1\cap S_2=\emptyset\}$,
where we use the shorthand
$S_1\uplus S_2\assign (S_1{\times}\{1\})\cup (S_2{\times}\{2\})$.

\heading{$\Z_2$-spaces and $\Z_2$-index. }
A \emph{$\Z_2$-space} (also called \emph{antipodality
space} in the literature) is a pair $(T,\nu)$, where $T$ is a topological
space and $\nu\:T\to T$, called the \emph{$\Z_2$-action},
 is a homeomorphism such that
$\nu^2=\nu\circ\nu=\id_T$. If $(T_1,\nu_1)$ and $(T_2,\nu_2)$
are $\Z_2$-spaces, a $\Z_2$-map between them is
a continuous mapping $f\:T_1\to T_2$ such that $f\circ\nu_1=\nu_2\circ f$.
The sphere $S^n$ is considered as a $\Z_2$-space with the
antipodal homeomorphism $x\mapsto -x$. Following
\v{Z}ivaljevi\'c \cite{Zivaljevic-topmeth}, we define the
 \emph{$\Z_2$-index} of
a $\Z_2$-space $(T,\nu)$ by
$$
\ind(T,\nu)\ \assign\ \min \left\{n\geq 0\sep \mbox{there is
a $\Z_2$-map $(T,\nu)\to S^n$}\right\}\in \{0,1,2,\ldots\}\cup\{\infty\}
$$
(the $\Z_2$-action $\nu$ is omitted from the notation if it is clear
from context).
If $\ind(T_1,\nu_1)>\ind(T_2,\nu_2)$, then there is no
$\Z_2$-map $T_1\to T_2$.
The Borsuk--Ulam theorem can be re-stated as \  $\ind(S^n)=n$.

A \emph{simplicial $\Z_2$-complex} is a 
simplicial complex $\K$
with a simplicial map $\nu$ of $\K$ into itself such that
(the canonical affine extension of) $\nu$ is a $\Z_2$-action
on $\polyh {\K}$. For the deleted join $\K^{*2}_\Delta$,
we have the canonical $\Z_2$-action given by ``swapping the two copies
of $V(\K)$,'' formally $(v,1)\mapsto (v,2)$ and $(v,2)\mapsto (v,1)$.

For any simplicial $\Z_2$-complex $\K$ whose $\Z_2$-action
is free (that is, has no fixed point), we have
\[
\dim\K\geq \ind\K \geq 1+\acyc( \K) \geq 1+\conn(\K).
\]
Here the first inequality needs freeness (in fact, $\ind\K=\infty$
if the $\Z_2$-action has a fixed point).
The second inequality is a homological version of the
Borsuk--Ulam theorem; see Walker \cite{Walker-bu}.
The parameter $\conn(\K)$ denotes the smallest $k$ such that
there exists a continuous map $S^{k+1}\to \polyh\K$ that 
is not nullhomotopic, while the acyclicity parameter is defined by
\[
\acyc(\K)\ \assign\  
\max\{k\sep \widetilde{\mathrm H}_i(\K,\Z_2)=0\textrm{ for all }i\le k\}.
\]
The $\Z_2$-acyclicity  is of interest in this context, since
it is effectively computable,
both theoretically  \cite[\S11]{Munkres} and practically
(for not too large complexes; see \cite{dumas:_comput_smith}),
while the $\Z_2$-index and the connectivity are in general
harder to determine.

We recall that two topological spaces $X$ and $Y$
are \emph{homotopy equivalent} if there are continuous maps
$f\:X\to Y$ and $g\:Y\to X$ such that $f\circ g$ is homotopic to
$\id_Y$ and $g\circ f$ is homotopic to $\id_X$. For $\Z_2$-spaces,
$\Z_2$-homotopy equivalence is defined analogously,
but we require that $f$, $g$, as well as all maps in the two homotopies
be $\Z_2$-maps.

\section{Proof methods for Kneser's conjecture}\label{s:overw}

\heading{The box complex $\Boxc G$. } For a graph $G$ and
any subset $A\subseteq V(G)$, let 
$$\CN(A)\ \assign\ \{v\in V(G)\sep
\{a,v\}\in E(G)\mbox{ for all } a\in A\}\ \ \subseteq\ \ V\setminus A
$$ 
be the set of all \emph{common neighbors} of $A$.

We define the \emph{box complex} $\Boxc{G}$ of a graph $G$
as the simplicial complex with vertex set $X=V(G){\times} [2]$
(\ie, two disjoint copies of $V(G)$), with simplices given by
\begin{eqnarray*}
\Boxc G&\assign&\bigl\{A'\uplus A''\sep
A',A''\subseteq V(G),\ A'\cap A''=\emptyset,\\
&&\hspace{21mm}G[A',A'']\mbox{ is complete}, \ 
\CN(A'),\CN(A'')\neq\emptyset\bigr\}.
\end{eqnarray*}
(We recall the notation 
$A'\uplus A''=(A'{\times}\{1\})\cup(A''{\times}\{2\})$.)
So the simplices of $\Boxc G$ correspond to complete bipartite subgraphs
in $G$. We admit $A'$ or $A''$ empty,
but then it is required that all vertices of the other shore
have a common neighbor (if both $A'$ and $A''$ are nonempty,
the condition $\CN(A'),\CN(A'')\neq\emptyset$ is superfluous).

However, if the extra condition on ``having a common neighbor''
is deleted, then we get a \emph{different} box complex
\[
\Bzero G\ \assign\ \bigl\{A'\uplus A''\sep
A',A''\subseteq V(G),\ A'\cap A''=\emptyset,\ G[A',A'']\mbox{ is
  complete}\bigr\}
\]
that contains $\Boxc G$, and which will also play a 
role in the following.

A canonical simplicial $\Z_2$-action on $\Boxc G$ is given by
interchanging the two copies of $V(G)$; that is,
$(v,1)\mapsto (v,2)$ and $(v,2)\mapsto (v,1)$, for $v\in V(G)$.
This makes $\Boxc G$ into a $\Z_2$-space.

If $f\:V(G)\to V(G)$ is a graph homomorphism, we associate to
it a map $\Boxc f\: V(\Boxc G)\to V(\Boxc H)$ in the obvious way:
$\Boxc f(v,j)\assign (f(v),j)$ for $v\in V(G)$, $j\in [2]$.
It is easily verified that $\Boxc f$ is a simplicial
$\Z_2$-map of $\Boxc G$ into $\Boxc H$. Moreover, the construction
respects the composition of maps, and so $\Boxc .$ can be
regarded as a functor from the category of graphs with
homomorphisms into the category of $\Z_2$-spaces with
$\Z_2$-maps. 

It is not hard to show that $\Boxc{K_m}$ is $\Z_2$-homotopy 
equivalent to $S^{m-2}$,
and  $\ind \Boxc{K_m} =m-2$ (see Section~\ref{sec:box_complexes}).
Since an $m$-coloring of $G$ can be regarded
as a homomorphism of $G$ into $K_m$, it induces a
$\Z_2$-map of $\Boxc G$ into $S^{m-2}$,
and we obtain
\begin{equation}\label{e:boxc-chi}
\chi(G)\geq \ind\Boxc G +2.
\end{equation}

The box complex $\Boxc G$ is a variation of ideas 
from
Alon, Frankl, and Lov\'asz \cite{AlonFranklLovasz} and \kriz~\cite{Kriz92}.

\heading{Neighborhood complexes and the  Lov\'asz bound. } 
Lov\'asz \cite{Lovasz-kneser}
defined the \emph{neighborhood complex} as 
$\N(G)\assign \{S\subseteq V(G)\sep \CN(S)\neq\emptyset\}$, and 
he proved that one always has
$\chi(G)\ge3+\conn(\N(G))$. His proof uses another
simplicial complex $\Lov{G}$, which can be defined as the
order complex of the system of all ``closed sets'' in~$\N(G)$:
$$
\Lov{G}\ \assign\ \Delta\{A\subset V(G)\sep \CN(\CN(A))=A\}.
$$
Thus, the vertices of $\Lov{G}$ are shores of inclusion-maximal
complete bipartite subgraphs of~$G$.
Unlike $\N(G)$, this $\Lov{G}$ is a simplicial $\Z_2$-complex,
with the $\Z_2$-action given by $A\mapsto \CN(A)$,
and a slight modification of
Lov\'asz' proof actually yields the lower bound
$$
\chi(G)\geq \ind \Lov{G}+2.
$$
As shown in \cite{Lovasz-kneser}, $\Lov{G}$ is a strong deformation
retract of $\N(G)$. 
A version of Lov\'asz' bound, formulated in terms of the classifying
map of the $\Z_2$-bundle associated with $\Lov G$,
was published by Milgram and Zvengrowski \cite{MilZve}.
Their formulation can be used as a tool for bounding 
$\ind \Lov{G}$ from below.

In Section~\ref{sec:box_complexes}
we will show that $\ind \Lov{G}=\ind \Boxc G$.
So while $\Boxc G$ and $\Lov{G}$ provide the same lower bound,
the functoriality of $\Boxc . $ (which was probably known
to experts, but as far as we know, hasn't appeared in print)
is a significant advantage.
Walker \cite{walker-kne} shows how a homomorphism induces
a $\Z_2$-map for the $\Lov.$ complexes, but the construction
is more complicated and not canonical.

\bigskip
The subsequent lower bounds are formulated for Kneser graphs,
in terms of the defining set system $\FF$.
The next definition is crucial in their formulation.
With a set system $\FF$, we associate
the following simplicial complex $\K=\K(\FF)$:
The vertex set of $\K$ is $X$, the ground set of $\FF$,
and
$$
\K(\FF)\ \assign\ \{S\subseteq X\sep
F\not\subseteq S\textrm{ for all } F\in\FF\}.
$$
Thus $\FF$ is the family of ``minimal nonfaces'' of~$\K$
(plus possibly additional nonfaces), 
while $\K$ is the complex of ``$\FF$-free sets.''

\heading{The Sarkaria bound. } From Sarkaria's
proof of Kneser's conjecture \cite{Sarkaria-kneser},
the following general bound can be deduced
($\FF$ is assumed to have the ground set $[n]$):
$$
\chi(\KG(\FF))\geq \ind\Delta\Bigl((2^{[n]})_\Delta^{*2}\setminus\K^{*2}_\Delta
\Bigr)+1.
$$
As it turns out, the complex on the right-hand side is just another
version of a box complex of $\KG(\FF)$. Sarkaria, in the concrete
cases he deals with, then proceeds to estimate the index 
of that complex using an elegant trick with joins
(``Sarkaria's inequality''; see \cite{Zivaljevic-topmeth} or
\cite{Matousek:BU}), which in general leads to
\begin{equation}\label{e:sark-b}
\chi(\KG(\FF))\geq n-1 - \ind \K^{*2}_\Delta.
\end{equation}
We call the right-hand side of this inequality the Sarkaria bound.
It is not explicitly stated in this way in Sarkaria's papers,
and so perhaps ``generalized Sarkaria bound'' would be more
precise, but repeating the adjective ``generalized'' at every
occasion seems annoying.

\heading{B\'ar\'any's proof } from \cite{Barany-kneser} yields a lower
bound that can generally be phrased as follows.
\emph{Suppose that for some $d\geq 1$,
the ground set $X$ of $\FF$ can be placed into the sphere
 $S^{d}$ in such a way that for every open hemisphere $H$
there exists a set $F\in\FF$ with $F\subseteq X\cap H$.
Then $\chi(\KG(\FF))\geq d{+}2$. }
(Kneser's conjecture is obtained from this using Gale's lemma,
stating that, for every $d,k\geq 1$, one can place
 $2k{+}d$ points on $S^d$ so that
every open hemisphere contains at least $k$ points.)
For the purposes of comparing the
 bound with the other
bounds, we will rephrase it using the Gale transform;
see Section~\ref{s:bar-gale}. The result can be expressed
as follows: \emph{Suppose that $\K$ is a subcomplex of the
boundary complex of an $(n{-}d)$-dimensional convex
polytope $P$ (under a suitable identification of the vertices of
$\K$ with the vertices of $P$).
Then $\chi(\KG(\FF))\geq d$. } We will refer to a number $d$ as in
this statement
as the B\'ar\'any bound; a comment similar to the one for 
the Sarkaria bound applies here as well.

{}From this form it is not hard to show that the
Sarkaria bound is always at least
as strong as the B\'ar\'any bound (but, of course, the index
in (\ref{e:sark-b}) might be difficult to evaluate).

\heading{The \dolnikov--\kriz\ bound } 
is a purely combinatorial lower estimate for $\chi(\KG(\FF))$.
For a set system $\FF$, let the
\emph{$2$-colorability defect} $\cd_2(\FF)$ (called the \emph{width}
in \cite{Kriz92})  be the minimum size of a subset $Y\subseteq X$
such that the system of the sets of $\FF$ that contain no points
of $Y$ is $2$-colorable. In other words,
we want to color each point of $X$
red, blue, or white in such a way that no set of $\FF$ is completely
red or completely blue (it may be completely white), and
$\cd_2(\FF)$ is the minimum number of white points required
for such a coloring. The following bound was derived
by \dolnikov\ \cite{dolnikov82,dolnikov88}
by a geometric argument from the Borsuk--Ulam theorem, and
independently (and as a part of a more general result)
by \kriz\ \cite{Kriz92,Kriz92err}, via certain box complexes:
\begin{equation}\label{e:do-kri}
\chi(\KG(\FF))\geq \cd_2(\FF).
\end{equation}
(Since it is easily seen that  
$\cd_2(\binom{[n]}{k})=n-2k+2$, Kneser's conjecture follows.)
A very short and elegant geometric reduction to a suitable
version of the Borsuk--Ulam theorem follows immediately
from the recent work of Greene \cite{Greene-kne}, which currently provides
the shortest self-contained proof of the Kneser conjecture.

The inequality (\ref{e:do-kri}) is also
an immediate consequence of (\ref{e:sark-b}).
Indeed, estimating the $\Z_2$-index by the dimension,
(\ref{e:sark-b}) leads to $\chi(G)\geq 
 n-1 - \dim \K^{*2}_\Delta$, and
some unwrapping of definitions
reveals that, surprisingly, the latter quantity is exactly~$\cd_2(\FF)$.

\section{The hierarchy}
In the following theorem, we summarize and compare all the considered
lower bounds for~$\chi(G)$. 

\begin{theorem}[The Hierarchy Theorem]\label{thm:Main}
Let $G=(V,E)=\KG(\FF)$ be a finite (Kneser) graph
with no isolated vertices,
where $\FF\subseteq 2^{[n]}$, and let 
$\K=\K(\FF)=\{S\subseteq [n]\sep F\not\subseteq S$ 
$\textrm{for all } F\in\FF\}$. Then 
we have the following chain\renewcommand\thefootnote{{\bfseries*}}%
\footnote{We have labelled the equations and inequalities in the 
following chain of by (\ref{ineq:1})--(\ref{ineq:12}); we will refer
to these labels below when we prove the relations, one by one.
Note that this is a chain of inequalities and equations, except at
the end, where we do not imply a relation between the B\'ar\'any bound
and the \dolnikov--\kriz\ bound.}
of inequalities and equalities:

\vbox{%
$$
{\renewcommand\arraystretch{1.3}
\begin{array}{rcll}
\displaystyle
\chi(G)&\mathop\ge\eqntag{ineq:1}&\ind \Boxc{G} +2 
       \ \mathop=\eqntag{ineq:4}\ \ind\Lov{G}+2 
                & {\textrm{``the Lov\'asz bound''}}    \\[4pt]
       &\mathop\ge\eqntag{ineq:8}&
        \ind \Bzero{G} +1\ \mathop=\eqntag{ineq:10}\ \ind\Delta\Bigl((2^{[n]})_\Delta^{*2}
                       \setminus\K^{*2}_\Delta\Bigr)+1 \\[4pt]
       &\mathop\ge\eqntag{ineq:11}& 
         n-1 - \ind \K^{*2}_\Delta &
        {\textrm{``the  Sarkaria bound''}}             \\[4pt]
       &&\mathop\ge\eqntag{ineq:bar}\quad d
           \quad \textrm{if $\K\subseteq\partial P$ for 
                                    an $(n{-}d)$-polytope $P$}
       &\textrm{``the B\'ar\'any bound''}             \\[2pt]
       &&\mathop\ge\eqntag{ineq:12}\quad 
        n-1 - \dim \K^{*2}_\Delta\ =\ \cd_{2}(\FF)&
         {\textrm{``the \dolnikov--\kriz\ bound.''}}  
\end{array}
}
$$
\begin{picture}(0,0)
\put(48,31){\Bigg\{}
\end{picture}}%
\end{theorem}

We have already proved (\ref{ineq:1}) 
(which is identical to (\ref{e:boxc-chi})),
as well as~(\ref{ineq:12}). The inequality (\ref{ineq:11}) 
was essentially proved by
Sarkaria (\cite{Matousek:BU} contains a detailed proof).
The B\'ar\'any inequality (\ref{ineq:bar}) is proved in
Section~\ref{s:bar-gale}, and the remaining
claims (\ref{ineq:4}), (\ref{ineq:8}) and (\ref{ineq:10}) follow from our
discussion of box complexes in Section~\ref{sec:box_complexes}.

\heading{Remarks on gaps/tightness. }
\begin{description}
\item[{\rm (\ref{ineq:1})}]
The gap in the first inequality may be
arbitrarily large. For example, for graphs without a $4$-cycle,
which can have arbitrarily large chromatic number,
all the topological lower bounds presented here are trivial.
Indeed, for every~$G$ without a $4$-cycle, there is a canonical
$\Z_2$-equivariant $\Z_2$-map of $\sd\Boxc{G}$ to a
$1$-dimensional subcomplex of $\Boxc G$, given by 
$A'\uplus A''\mapsto\emptyset\uplus\CN(A')  $ for $|A' |\ge2$,
$A'\uplus A''\mapsto\CN(A'') \uplus\emptyset$ for $|A''|\ge2$,
and $A'\uplus A''\mapsto A'\uplus A''$ otherwise, 
and we get $\ind\Boxc{G}\le1$.
\item[{\rm (\ref{ineq:8})}]
The gap in the inequality {\rm (\ref{ineq:8})} can be at most~$1$;
this can be derived from the inequality (\ref{m:i}) of
Proposition~\ref{prop:box_complexes} below, which yields
$\ind\Bzero{G}\ge\ind\Boxc{G}$.
\item[{\rm (\ref{ineq:bar})}]
The   B\'ar\'any bound can be strictly larger
than the \dolnikov--\kriz\ bound, 
without any bound on the gap, as will
be discussed at the end of Section~\ref{s:bar-gale}
for the example of the Schrijver graphs.
Thus, in particular, the gap in~(\ref{ineq:8}) can be arbitrarily large.

The B\'ar\'any bound depends on the choice of the polytope $P$.
At present we do not know whether $P$ can always be chosen of
dimension at most $\dim \K^{*2}_\Delta+1$, that is, 
whether the B\'ar\'any bound can always be made at least as
strong as the \dolnikov--\kriz\ bound. 
On the other hand, we do not have an example
where the B\'ar\'any bound is necessarily smaller
than the Sarkaria bound. 
\end{description}

\heading{Remarks on size/computability. }
Although the \dolnikov-\kriz\ ``colorability defect'' bound
is attractive since it is \emph{combinatorial},
in general the Lov\'asz bound may be much tighter. However,
the number of vertices of $\Lov G$ may be exponential in $n=|V(G)|$,
and similarly for some of the other box complexes. On the other hand,
$\Boxc G$ has only $2n$ vertices. The number of simplices can still be
exponential, but if, for example, the maximum degree of $G$ is bounded
by a constant, then there are at most polynomially many simplices.
Perhaps a  computation of the $\Z_2$-acyclicity of $\Boxc G$,
which provides a lower bound for $\ind\Boxc G$ (and thus for $\chi(G)$),
might be feasible in some cases.

\section{Box complexes and neighborhood complexes}\label{sec:box_complexes}

In the following definition, we collect six (natural)
variants of box complexes,
four defined for a graph and two for a Kneser representation
of it. For completeness, we also include the
box complexes $\Boxc G$ and $\Bzero G$ that were already defined above.

\begin{definition}[Box complexes]
Let $G=(V,E)=\KG(\FF)$ be a finite (Kneser) graph with no
isolated vertices, and suppose that the ground set of $\FF$ is
$[n]$.
The first two complexes are on the vertex set $V{{\times}}[2]$;
they were already defined in Section~\ref{s:overw}.
\begin{enumerate}
\item
The box complex $\Boxc G$ is
\begin{eqnarray*}
\Boxc G&\assign&\bigl\{A'\uplus A''\sep
A',A''\subseteq V,\ A'\cap A''=\emptyset,\\
&&\hspace{21mm} G[A',A'']\mbox{ is complete},
\ \CN(A'),\CN(A'')\neq\emptyset\bigr\}.
\end{eqnarray*}
Equivalently, but more concisely, we can also write
$\Boxc G=\bigl\{A'\uplus A''\sep A' \subseteq \CN(A'')\neq\emptyset$,
                    $A''\subseteq \CN(A' )\neq\emptyset\bigr\}$.

\item
A simpler definition, but a larger complex, is obtained as 
\[
\Bzero G\ \assign\ \bigl\{A'\uplus A''\sep A',A''\subseteq V,\
A'\cap A''=\emptyset,\ G[A',A'']\mbox{ is complete}\bigr\}.
\]
This is almost as for $\Boxc G$, but here if one shore is empty, the other
can be anything.

\item
The following definition of a box complex,
from \kriz\ \cite[p.~568]{Kriz92}, takes into account only
the complete bipartite graphs with \emph{both} shores $A'$ and $A''$ nonempty:
\begin{eqnarray*}
\Bone G&\assign&
\Delta\bigl\{A'\uplus A''\sep \emptyset\neq A',A''\subset V, \\
&&\hspace{24mm} A'\cap A''=\emptyset,\ G[A',A'']\mbox{ is complete}\bigr\}.
\end{eqnarray*}
Here the \emph{vertices} are the vertex sets of complete bipartite
subgraphs of $G$, and the simplices are chains of such sets under
inclusion.
\item The vertices of the next box complex,
from Alon, Frankl \& Lov\'asz \cite[p.~361]{AlonFranklLovasz},
are directed edges
of $G$; that is, ordered pairs $(u,v)$ with $\{u,v\}\in E$.
We let
\begin{eqnarray*}
\Bedge{G}&\assign& \bigl\{\vec {F}\subseteq A'{\times} A'' \sep
\emptyset\neq A',A''\subset V,\ A'\cap A''=\emptyset,\ 
 G[A',A'']\mbox{ is complete}\bigr\}.
\end{eqnarray*}
That is, simplices are subsets of edge sets of complete bipartite
subgraphs of $G$, where the edges are oriented from the first
shore to the second shore. 
\item The simplicial complex 
$\Delta\bigl((2^{[n]})_\Delta^{*2}\setminus\K^{*2}_\Delta\bigr)$
appearing in Theorem~\ref{thm:Main}, used as an intermediate step
by Sarkaria in the derivation of his lower bounds,
can be more explicitly written as
\begin{eqnarray*}
\Bsark{\FF}&\assign&\Delta\bigl\{B'\uplus B''\sep
B',B''\subseteq [n],\ B'\cap B''=\emptyset, \\
&&\hspace{24mm} \mbox{ at least one of $B',B''$
contains a set of $\FF$}\bigr\}.
\end{eqnarray*}
The vertex set are pairs of disjoint subsets of the
ground set of $\FF$  that support a complete bipartite
subgraph of the Kneser graph, with at least one shore nonempty.
\item Finally, another Kneser box complex,
as in  \kriz\ \cite[p.~574]{Kriz92}, is
\begin{eqnarray*}
\Bonekneser{\FF}&\assign& \Delta\bigl\{B'\uplus B''\sep
B',B''\subseteq [n],\ B'\cap B''=\emptyset, \\
&&\hspace{24mm} \mbox{ both $B',B''$ contain a set of $\FF$}\bigr\}.
\end{eqnarray*}
\end{enumerate}
\end{definition}

On each of these types of box complexes, we have the
natural $\Z_2$-action that interchanges the shores of the bipartite subgraph.

As we will show, all these box complexes fall into two groups,
and those in each group have the same $\Z_2$-index.
Moreover, the Lov\'asz complex $\Lov{G}$ can also be included
in one of the groups.

\begin{theorem}
The following holds for the $\Z_2$-indices of the various box complexes:
\begin{eqnarray*}
       \ind\Bone G
& =&   \ind\Bonekneser\FF
\ =\   \ind\Bedge G 
\ =\   \ind\Boxc G
\ =\   \ind\Lov G
\\[4pt]
& \le& \ind\Bzero G
\ =\   \ind\Bsark\FF
\ \le\ \ind\Bone G +1.
\end{eqnarray*}
\end{theorem}

For a proof of this theorem, the following proposition provides
explicit simplicial $\Z_2$-maps among the various box complexes.
Here $\susp\K$ denotes the \emph{suspension} of a simplicial
complex $\K$ (a ``double cone'' over $\K$):
$\susp\K\assign\K\cup\{S\cup\{s\}\sep S\in\K\}\cup
\{S\cup\{n\}\sep S\in\K\}$, where $s$ and $n$ are two new
vertices not belonging to $V(\K)$.
The last inequality of the theorem follows from the
map (\ref{m:iv}) in the proposition,
together with $\ind\susp \K\leq\ind\K+1$
and the fact that $\susp S^n$ is an $S^{n+1}$.

\begin{proposition}[$\Z_2$-maps]\label{prop:box_complexes}
For every finite graph $G=\KG(\FF)$ without isolated vertices,
there are canonical simplicial $\Z_2$-maps
{%
\newcounter{saveeqn}%
\setcounter{saveeqn}{\theequation}%
\def\theequation{M\arabic{equation}}
\setcounter{equation}{0}
\begin{eqnarray}
\Boxc G &\longrightarrow & \Bzero G,\label{m:i}\\
\sd\Bedge{G}&\longleftrightarrow& \Bone {G},\label{m:ii}\\
\Bone{G}&\longrightarrow &\sd\Boxc{G},\label{m:iii}\\
\sd\Bzero{G}&\longrightarrow&\susp\Bone{G},\label{m:iv}\\
\Bone{G}&\longleftrightarrow&\Bonekneser{\FF},\label{m:v}\\
\sd\Bzero{G}&\longleftrightarrow&\Bsark{\FF},\label{m:vi}\\
\sd\sd \Boxc G&\longrightarrow &\Bone G,\label{m:vii}\\
\sd\Lov G&\longrightarrow &\Bone G,\label{m:viii}\\
\sd\Bone G&\longrightarrow &\sd\Lov G.\label{m:ix}
\end{eqnarray}
\setcounter{equation}{\thesaveeqn}
}\vskip-4mm
\end{proposition}

Answering a question by the first author,
Lov\'asz (personal communication, February 2000)
proved that $\Lov G$ is homotopy equivalent to $\Bone G$
(using the nerve theorem);
our construction of the $\Z_2$-map in (\ref{m:ix})
is inspired by his proof.

After a preliminary version of this paper was written,
Csorba \cite{Csorba-susp} proved that $\Bzero G$ is
$\Z_2$-homotopy equivalent to $\susp(\Boxc G)$ for every $G$,
thereby considerably strengthening our statements above concerning these
box complexes.

The following figure sketches the three main box complexes
for $G=C_5$ and suggests the maps between them. Here $\Boxc{C_5}$ is
homeomorphic to $S^1\times I$;
it is a subcomplex of $\Bzero{C_5}$, which additionally contains
two simplices on $5$ vertices. The complex $\Bone{C_5}$, an
$S^1$ on $10$ vertices, embeds into the barycentric
subdivision of~$\Boxc{C_5}$:

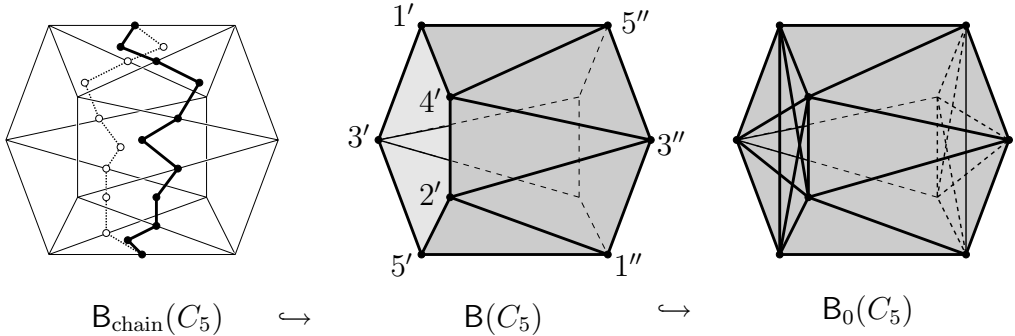
\begin{center}
\input{boxcomplexes.pstex_t}
\end{center}

\proof
\begin{enumerate}
\item[(\ref{m:i})]
This map is simply the identity on the vertex set: an inclusion map.
\item[(\ref{m:ii})]
The map $\sd\Bedge{G}\rightarrow\Bone {G}$ 
is defined on the vertices of $\sd\Bedge{G}$
by setting $\vec F\mapsto h(\vec F)\uplus
t(\vec F)$, where 
$h(\vec F)$ collects the tails, and
$t(\vec F)$ collects the heads, of the directed edges in~$\vec F$. 
 
The map $\Bone {G}\rightarrow\sd\Bedge{G}$
in the other direction is obtained by 
sending $A'\uplus A''$ to the (directed) edge set
of the complete bipartite graph with shores $A'$ and $A''$.

\item[(\ref{m:iii})]
The map $\Bone G\rightarrow\sd\Boxc G$ is 
simply an inclusion, as each vertex of $\Bone G$ is also a
simplex of $\Boxc G$.

\item[(\ref{m:iv})]
The map $\sd\Bzero{G}\rightarrow\susp\Bone{G}$
is obtained by mapping the ``improper'' complete bipartite
subgraphs (with one shore empty) to the suspension points.

\item[(\ref{m:v})]
A canonical $\Z_2$-map $\Bonekneser{\FF}\rightarrow\Bone{G}$ 
is defined on the vertices by
\[
B'\uplus B''\ \longmapsto\ \{F'\in\FF\sep
F'\subseteq B'\}\uplus 
\{F''\in\FF\sep F''\subseteq B''\}.
\]
A canonical $\Z_2$-map $\Bone{G}\rightarrow\Bonekneser{\FF}$ is 
obtained by mapping the vertices: 
\[
A'\uplus A''\ \longmapsto\ \textstyle\big(\bigcup A'\big)\uplus 
\big(\bigcup A''\big).
\]
\item[(\ref{m:vi})]
The same formulas define  $\Z_2$-maps
$\Bsark{\FF}\longleftrightarrow \sd\Bzero{\KG(\FF)}$.
\item[(\ref{m:vii})]
The  map $\sd\sd \Boxc G\to\Bone G$
 is defined on the second barycentric
subdivision of $\Boxc G$, and so its \emph{vertices} are chains
of the form 
$$
\AA=
\bigl(A'_0\uplus A''_0\subset\cdots\subset \cdots\subset A_k'\uplus A_k''\bigr).
$$
We let $\mu'(\AA)$ be the smallest nonempty set in the chain of sets
$\AA'\assign (A_0'  \subset\cdots\subset   A_k'  \subseteq
\CN(A_k'')\subseteq\cdots\subseteq \CN(A_0''))$,
and similarly for $\mu''(\AA)$ (here we use the condition
$\CN(A'),\CN(A'')\neq\emptyset$ from the definition of $\Boxc G$).
We let the image of $\AA$ be $\mu'(\AA)\uplus\mu''(\AA)$.
This is a vertex of $\Bone{G}$: Since the first barycentric subdivision
does not contain $\emptyset$ as a vertex, at least
one of $A_0',A_0''$, say $A_0'$, is nonempty.
Then we have $\mu'(\AA)=A_0'$, while $\mu''(\AA)$
is contained in~$\CN(A_0'')$.
If we extend the
chain $\AA$, then this also leads to an extension of
the chains $\AA'$ and $\AA''$,
so the sets $\mu'(\AA)$ and $\mu''(\AA)$ can only
get smaller. Therefore, our map is simplicial, and
it is clearly a $\Z_2$-map.
\item[(\ref{m:viii})]
We recall that the vertices of
$\Lov{G}$ are the nonempty subsets $A\subset V$ that are \emph{closed}
in the sense that $A=\CN(\CN(A))$, or equivalently, $A=\CN(B)$ for some
nonempty subset $B\subset V$.

First we define a $\Z_2$-map $f\:\sd \Lov{G}\to\Bone G$.
A vertex of $\sd\Lov{G}$ is a chain
$\AA=(A_0\subset A_1\subset\cdots\subset A_k)$
of nonempty closed sets. We set $f(\AA)\assign A_0\uplus \CN(A_k)$.
Since $\CN(A_k)\subseteq \CN(A_0)$,
the image is indeed a vertex of $\Bone G$.
If a chain $\AA'$ extends $\AA$, its first set can only be smaller
than the first set of $\AA$, and the last set can only be larger
than the last set of $\AA$. Therefore, $f(\AA')\subseteq
f(\AA)$, and it follows that $f$ is simplicial.
Finally, the image of $\AA$ under the $\Z_2$-action on $\sd\Lov{G}$
is the chain $\BB=(\CN(A_k)\subset \CN(A_{k-1})\subset\cdots
\subset\CN(A_0))$. We have $f(\AA)=A_0\uplus \CN(A_k)$
and $f(\BB)=\CN(A_k)\uplus\CN^2(A_0)=\CN(A_k)\uplus A_0$
(as $A_0$ is closed), and so $f$ is a $\Z_2$-map.

\item[(\ref{m:ix})]
Finally, we provide a $\Z_2$-map 
$\sd\Bone G\to\sd\Lov G$.
A vertex in $\sd\Bone G$ is a chain
$\AA=\big(A_0'\uplus A_0''\subset\cdots\subset A_k'\uplus A_k''\big)$.
All the sets $\CN^2(A'_0)$,\ldots, $\CN^2(A'_k)$,
$\CN(A''_0)$,\ldots, $\CN(A''_k)$ are closed and nonempty,
and the following inclusion are easily verified:
$\CN^2(A'_0)\subseteq \cdots\subseteq \CN^2(A'_k)\subseteq
\CN(A''_k)\subseteq\cdots\subseteq \CN(A''_0)$.
So by omitting repeated sets from this chain,
we obtain a vertex of $\sd\Lov G$.
The chain $\AA$ is mapped to this vertex.
If we extend $\AA$, the image stays the same or is extended
as well, so the map is simplicial. Finally, it is a
$\Z_2$-map; here we use that $\CN^3=\CN$.
\proofend
\end{enumerate}

The above definitions by far do not exhaust the list of 
possible (and possibly interesting) box complex variants. 
For example, we could also define the
``Kneser counterparts'' of $\Boxc G$ and $\Bzero G$, namely 
$\Boxckneser{\FF}\assign\bigl\{S\subseteq B'\uplus B''\sep
B',B''\subseteq [n], B'\cap B''=\emptyset$,
both $B'$ and $B''$ contain a set of $\FF\bigr\}$
and
$\Bzerokneser{\FF}\assign\bigl\{S\subseteq B'\uplus B''\sep
B',B''\subseteq [n], B'\cap B''=\emptyset$,
at least one of $B',B''$ contains a set of $\FF\bigr\}$.
By similar arguments as in the proof of Proposition~\ref{prop:box_complexes},
it can be shown that $\ind\Bzero{\KG(\FF)}=\ind\Bzerokneser \FF$
and $\ind\Boxc{\KG(\FF)}=\ind\Boxckneser \FF$,
but we prefer to omit this part.

In some of the cases in Proposition~\ref{prop:box_complexes},
the maps even provide $\Z_2$-homotopy equivalences between the respective
complexes. Since the proofs are not very interesting and, at present, 
the  $\Z_2$-homotopy types do not seem to bring anything
new concerning the lower bounds for the chromatic number, we have decided
not to discuss this in the present paper.

\heading{The box complexes of complete graphs. }
In order to use a box complex for bounding the chromatic number
of a graph, we need to know the $\Z_2$-index of the box complex
of $K_m$. In view of the equalities of indices established above
for the various box complexes, we mention only two of the box
complexes here (but the others can also be analyzed directly
without difficulty).

\begin{proposition}\label{prop:color_spaces}
The polyhedron of the simplicial complex $\Bzero{K_m}$ is 
(homeomorphic to) an $S^{m-1}$, and thus \ $\ind\Bzero{K_m}=m{-}1$.

The polyhedron of the simplicial complex $\Boxc{K_m}$
is homeomorphic to $S^{m-2}{\times} [0,1]$, and thus \ $\ind\Boxc{K_m}=m-2$.
\end{proposition}

\proof
$\Bzero{K_m}$ is isomorphic to the boundary
complex of an $m$-dimensional cross polytope (``generalized
octahedron,'' or unit ball of the $\ell_1$-norm), with its
canonical antipodal $\Z_2$-action.

The complex $\Boxc{K_m}$  is isomorphic to the boundary
complex of an $m$-dimensional cross polytope with two opposite facets
removed, and with its canonical antipodal $\Z_2$-action.
\proofend

By now, all parts of Theorem~\ref{thm:Main} have been proved,
except for the claims involving the  B\'ar\'any bound. 

\section{The  B\'ar\'any bound}\label{s:bar-gale}

\heading{Proof of the inequality (\ref{ineq:bar}).} We  need to prove
that if
the simplicial complex $\K$ is a part of the boundary of 
an $k$-dimensional convex polytope $P$, then
$\ind \K^{*2}_\Delta\leq k{-}1$. 
Briefly speaking, the reason is that the deleted join of $S^{k-1}$
contains $\polyh{\K^{*2}_\Delta}$ and 
is $\Z_2$-homotopy equivalent to $S^{k-1}$.
A more detailed argument follows.

Let us recall that if $X$ and $Y$ are topological spaces,
the join $X*Y$ is the quotient of the product space 
$X{\times}Y{\times}[0,1]$ by the equivalence $\approx$,
where $(x,y,0)\approx (x',y,0)$ and $(x,y,1)\approx (x,y',1)$
for all $x,x'\in X$ and $y,y'\in Y$.
If $\K$ is a simplicial complex and 
$\K^{*2}\assign\{S\uplus S'\sep S,S'\in\K\}$, then
there is a (canonical) homeomorphism $\polyh{\K}^{*2}\cong\polyh{\K^{*2}}$.

Here we need the space
$$
Z\ \assign\ \big(S^{k-1}*S^{k-1}\big)
\setminus \{(x,x,\tfrac12)\sep x\in S^{k-1}\};
$$
this is a kind of deleted join of $S^{k-1}$ (considered as
a topological space, not a simplicial complex).
The space $Z$ is equipped with the $\Z_2$ action
given by $(x,y,t)\mapsto (y,x,1{-}t)$.

Using the canonical homeomorphism
 $\polyh{\K^{*2}}\to \polyh{\K}^{*2}$ and a homeomorphism
of the boundary of $P$ with $S^{k-1}$,
we obtain a $\Z_2$-map $h\:\polyh{\K^{*2}}\to S^{k-1}*S^{k-1}$.
Moreover, if we restrict the left-hand side to $\polyh{\K^{*2}_\Delta}$,
then the image contains no point of the form
$(x,x,\frac 12)$, and so $h$ is a 
$\Z_2$-map $\polyh{\K^{*2}_\Delta}\to Z$.
So it suffices to show that $\ind Z\leq k{-}1$.
Let  $S^{k-1}$ be represented as the unit sphere in $\R^k$.
We  define a $\Z_2$-map $Z\to S^{k-1}$ by
$$
(x,y,t)\mapsto \frac{tx-(1{-}t)y}{\|tx-(1{-}t)y\|}.
$$
If $tx-(1{-}t)y=0$ for unit vectors $x,y$, then
$t=\frac 12$ and $x=y$. Thus, the map is well defined,
and one can also check that it is continuous.
Clearly, it commutes with the respective $\Z_2$-actions.
This concludes the proof.
\proofend

\heading{An extension of B\'ar\'any's proof. }
Here we directly prove the inequality
$\chi(\KG(\FF))\geq d$ whenever $\K(\FF)$ is a part of the boundary
of an $(n{-}d)$-dimensional simplicial 
convex polytope. This explains the relation
of this bound to B\'ar\'any's \cite{Barany-kneser} original proof 
of the Kneser conjecture.

Let $V=(v_1,v_2,\ldots,v_n)
\subset \R^{n-d}$ be the vertex set of $P$, and let 
$V^*=(v_1^*,v_2^*,\ldots,v_n^*)\subset \R^{d-1}$ be a Gale diagram
of $V$ (see, \eg, \cite[Lect.~6]{Ziegler-poly} 
for an introduction to Gale diagrams).
Without loss of generality, we may assume that
all points in $V^*$ different from $0$ lie in $S^{d-2}$.
 We recall that the points of $V$, and thus
also the points of $V^*$, are in one-to-one correspondence
with the elements of the ground set of $\FF$. 

We want to show that \emph{every open hemisphere in this $S^{d-2}$
contains a set corresponding to a set of $\FF$}. 
Once this is done, we can proceed
exactly as in B\'ar\'any's proof. Namely, supposing that the sets
of $\FF$ have been colored by at most $d{-}1$ colors,
we define the set $A_i\subseteq S^{d-2}$ as the set of all
$x\in S^{d-2}$ such that the open hemisphere centered at $x$
contains at least one set of $\FF$ colored by $i$, $i\in [d{-}1]$.
Each $A_i$ is open, and by the above claim,
the $A_i$ together cover all of $S^{d-2}$.
By a suitable version of the Borsuk--Ulam theorem
(called the \lju--\snirel\  theorem),
there exist $x\in S^{d-2}$ and $i\in [d{-}1]$
such that $x,-x\in A_i$. This means that two opposite
open hemispheres both contain a set of color $i$,
and so the coloring is not a proper coloring of the Kneser
graph.

To prove the claim, consider an open hemisphere $H$,
and let $S^*\assign H\cap V^*$. Since $H$ is defined by
an open halfspace, there is a linear functional on 
$\R^{d-1}$ that is positive on $S^*$ and nonpositive on
$V^*\setminus S^*$. 
Let $S\subset V$ be the set corresponding to $S^*$.
Properties of the Gale diagram imply
that there is an affine dependence of the points of $V$
in which the points of $S$ have positive coefficients
and the other points have nonpositive coefficients.
This means that $\conv(S)\cap\conv(V\setminus S)\neq\emptyset$.
So $S$ is not the vertex set of a face of $P$, and thus
$S\not\in\K(\FF)$. This means that $S$ contains a set of $\FF$.
The proof is finished.
\proofend

\heading{On Schrijver graphs and the \dolnikov--\kriz\ bound. }
Let $0<2k<n$, and let $\binom{[n]}k_\stab$ denote the system
of all sets $F\subseteq [n]$ such  that if $i\in F$ then
$i{+}1\not\in F$, and if $n\in F$ then $1\not\in F$.
(So the sets of $\binom{[n]}k_\stab$ can be identified with the independent
sets in the cycle of length $n$, with the numbering of vertices 
from $1$ to $n$ along the cycle.)
Schrijver \cite{Schrijver-kneser} proved that
the graph $\SG_{n,k}\assign\KG\bigl(\binom{[n]}k_\stab\bigr)$
is a vertex-critical subgraph of the Kneser graph
$\KG\bigl(\binom{[n]}k\bigr)$. That is,
$\chi(\SG_{n,k})=n-2k+2$, and every proper induced subgraph of $\SG_{n,k}$
has a smaller chromatic number. The criticality follows
by a clever coloring construction and we will not consider it
here; we look at a proof of the lower bound.

It turns out that the B\'ar\'any bound applies very neatly here.
Let $P\assign C_{2k-2}(n)$ be a cyclic polytope of dimension $2k-2$
on $n$ vertices. 
With the usual numbering of the vertices,  Gale's evenness
criterion (see \cite[Thm.~0.7]{Ziegler-poly}) shows that
the $S\subseteq [n]$ that contain no set of
$\binom{[n]}k_\stab$ are exactly the proper faces of~$P$.
Therefore, $\K\bigl(\binom{[n]}k_\stab\bigr)=\partial P$,
and the B\'ar\'any bound immediately yields $\chi(\SG_{n,k})\geq n-2k+2$.

The Schrijver graphs also provide examples where the
\dolnikov--\kriz\ bound is considerably
weaker than the B\'ar\'any bound. Indeed, it is easy to check that
$\cd_2\bigl(\binom{[n]}k_\stab\bigr)=n-4k+4$.

\section{Concluding remarks}

\setlength{\leftmargini}{5mm}
\begin{enumerate}
\item 
Many of the above considerations can easily be extended to 
Kneser hypergraphs (where vertices are again the sets of~$\FF$,
and edges are $r$-tuples of pairwise disjoint sets, for some
given $r$), and even to  the
$s$-disjoint Kneser hypergraph version of Sarkaria, as in 
Ziegler \cite{Z77}. A detailed exploration of this is a subject
for further research.

Indeed, the $p$-partite versions of some of our box complexes
as presented in Section~\ref{sec:box_complexes} appear in
the published literature that concerns Kneser hypergraphs:
The $p$-partite version of $\Bedge{G}$ is used in
Alon, Frankl \& Lov\'asz \cite[p.~361]{AlonFranklLovasz}.
The $p$-partite version of $\Bone{G}$ appears in \kriz\ \cite[p.~568]{Kriz92},
while in \kriz\ \cite[p.~574]{Kriz92} we find the $p$-partite
version of $\Bonekneser{G}$. The box complexes appear in
\kriz' work as ``resolution complexes'' for equivariant cohomology.

\item
In a discussion with Jarik Ne\v{s}et\v{r}il, we noted the following
interpretation of $\Bzero{G}$ (and $\Boxc{G}$):
the simplices of $\Bzero{G}$ are vertex sets of complete
bipartite subgraphs in $G\times K_2$, where the ``categorical product''
of graphs $G\times H$ has vertex set $V(G)\times V(H)$ and
edges $\{(u,v),(u',v')\}$ such that $\{u,u'\}\in E(G)$ and
$\{v,v'\}\in E(H)$. The $\Z_2$-action is then induced by the exchange
of vertices of $K_2$. A homomorphism $G\rightarrow G'$ induces a homomorphism
$G\times K_2 \rightarrow G' \times K_2$
(so the $\Boxc{\cdot}$ functor ``factors'' in this way).
This can be generalized, for example,
 to a product $G\times C_p$, where $C_p$ is the $p$-cycle with the
natural $\Z_p$-action (cyclic shift). 

\item
In a similar spirit but earlier,
Lov\'asz and others have considered a setting of ``Hom-complexes'':
For this let $H$ and $G$ be finite graphs, and
let $\Delta_G$ be a simplex with vertex set~$V(G)$.
The Cartesian power $(\Delta_G)^{|V(H)|}$ is a polyhedral complex,
whose vertices can be identified with mappings $V(H)\to V(G)$.
Then ${\rm Hom}(H,G)$ is the subcomplex induced by the vertices
that are homomorphisms (thus, the faces are all polyhedra 
$F\in(\Delta_G)^{|V(H)|}$ such that each vertex of $F$ is a homomorphism).
This rather general setting covers Lov\'asz' connectivity
lower bound, as well as some of the recent work of Brightwell and
Winkler~\cite{brightwell:_graph}. 
Lov\'asz conjectured that if ${\rm Hom}(C_{2r+1},G)$ is
$k$-connected then $\chi(G)\geq k+4$. Here $r\geq 1$ is an arbitrary
integer and $C_{2r+1}$ denotes the odd cycle of length $2r+1$.
This was recently proved by Babson and Kozlov \cite{BabKozl} by
advanced topological methods. Let us remark that while
$\chi(G)\geq k+3$ is easy to establish by methods discussed in
the present paper, the improvement by 1 currently seems hard,
and it represents new kind of topological obstruction to
$(k+4)$-colorability, which apparently is not captured by any
of the lower bounds discussed above.

\item
The inequality $\ind\susp\K \le \ind\K$ may
be strict for finite simplicial complexes $\K$:
For example, one may obtain a cell complex model by
taking $h\:S^3\to S^2$ to be the Hopf map, and
attaching two $4$-cells to $S^2$ via $2h$ resp.\ $-2h$,
where multiples of maps are taken according to
addition in $\pi_3(S^2)\cong\Z$.
(This particular example was suggested
by P\'eter Csorba, based on an earlier construction
by Csorba, \v{Z}ivaljevi\'c, and the first author
for a different purpose; another approach was
proposed by Wojchiech Chach\'olski.)

However, it is not clear whether $ \ind\susp\K =\ind \K$ occurs
in the rather special setting of~(\ref{ineq:8}), where $\K$ is a box
complex. An interesting step in this direction was recently taken
by Csorba (private communication), who proved that given any 
simplicial $\Z_2$-complex $\K$, there is a graph $G$
such that $\Nbhd G$ is homotopy equivalent to $\K$.
At present it is not clear whether this result can be extended
to a $\Z_2$-homotopy equivalence of the box complex of a suitable
graph with a given $\K$; if yes, this would provide
many pathological examples and, in particular, it would
answer the above question.

\end{enumerate}

\subsection*{Acknowledgements}
We would like to thank Imre B\'ar\'any,
Anders Bj\"orner,  L\'aszl\'o Lov\'asz,
and Jaroslav Ne\v{s}et\v{r}il for most valuable conversations
and comments.
We also thank Wojchiech Chach\'olski and P\'eter Csorba
for suggestions of examples with $\ind\K=\ind\susp\K$.

\end{document}

%% file: boxcomplexes.pstex_t
\begin{picture}(0,0)%
\includegraphics{boxcomplexes.pstex}%
\end{picture}%
\setlength{\unitlength}{2368sp}%
\begingroup\makeatletter\ifx\SetFigFont\undefined%
\gdef\SetFigFont#1#2#3#4#5{%
  \reset@font\fontsize{#1}{#2pt}%
  \fontfamily{#3}\fontseries{#4}\fontshape{#5}%
  \selectfont}%
\fi\endgroup%
\begin{picture}(10556,3475)(1039,-4796)
\put(9601,-4636){\makebox(0,0)[lb]{\smash{\SetFigFont{12}{14.4}{\familydefault}{\mddefault}{\updefault}{\color[rgb]{0,0,0}$\Bzero{C_5}$}%
}}}
\put(7876,-4636){\makebox(0,0)[lb]{\smash{\SetFigFont{12}{14.4}{\familydefault}{\mddefault}{\updefault}{\color[rgb]{0,0,0}$\hookrightarrow$}%
}}}
\put(3901,-4711){\makebox(0,0)[lb]{\smash{\SetFigFont{12}{14.4}{\familydefault}{\mddefault}{\updefault}{\color[rgb]{0,0,0}$\hookrightarrow$}%
}}}
\put(1951,-4711){\makebox(0,0)[lb]{\smash{\SetFigFont{12}{14.4}{\familydefault}{\mddefault}{\updefault}{\color[rgb]{0,0,0}$\Bone{C_5}$}%
}}}
\put(7426,-4186){\makebox(0,0)[lb]{\smash{\SetFigFont{12}{14.4}{\familydefault}{\mddefault}{\updefault}{\color[rgb]{0,0,0}$1''$}%
}}}
\put(7876,-2911){\makebox(0,0)[lb]{\smash{\SetFigFont{12}{14.4}{\familydefault}{\mddefault}{\updefault}{\color[rgb]{0,0,0}$3''$}%
}}}
\put(7501,-1636){\makebox(0,0)[lb]{\smash{\SetFigFont{12}{14.4}{\familydefault}{\mddefault}{\updefault}{\color[rgb]{0,0,0}$5''$}%
}}}
\put(5626,-3436){\makebox(0,0)[rb]{\smash{\SetFigFont{12}{14.4}{\familydefault}{\mddefault}{\updefault}{\color[rgb]{0,0,0}$2'$}%
}}}
\put(5326,-4186){\makebox(0,0)[rb]{\smash{\SetFigFont{12}{14.4}{\familydefault}{\mddefault}{\updefault}{\color[rgb]{0,0,0}$5'$}%
}}}
\put(4876,-2836){\makebox(0,0)[rb]{\smash{\SetFigFont{12}{14.4}{\familydefault}{\mddefault}{\updefault}{\color[rgb]{0,0,0}$3'$}%
}}}
\put(5101,-1561){\makebox(0,0)[lb]{\smash{\SetFigFont{12}{14.4}{\familydefault}{\mddefault}{\updefault}{\color[rgb]{0,0,0}$1'$}%
}}}
\put(5626,-2461){\makebox(0,0)[rb]{\smash{\SetFigFont{12}{14.4}{\familydefault}{\mddefault}{\updefault}{\color[rgb]{0,0,0}$4'$}%
}}}
\put(5851,-4711){\makebox(0,0)[lb]{\smash{\SetFigFont{12}{14.4}{\familydefault}{\mddefault}{\updefault}{\color[rgb]{0,0,0}$\Boxc{C_5}$}%
}}}
\end{picture}